\documentclass{article}
\usepackage{amscd}
\usepackage[T2A]{fontenc}
\usepackage[cp1251]{inputenc}
\usepackage[ukrainian]{babel}
\usepackage{amssymb} 
\newtheorem{theorem}{Theorem}[section]
\newtheorem{lemma}{Lemma}[section]

\newtheorem{cor}{Corralalry}[section]
\newtheorem{rem}{Remark}[section]
\newtheorem{defn}{Defintion}[section]

\begin{document}

\begin{center}
\textbf{\Large{ Quantization of Lyapunov functions}}\\
\end{center}

\begin{center}
 Yulia Sharko
\end{center}

We investigate discret conditions for stability and asymptotic stability by Lyapunov and the point of equilibrium of autonomous system of differential equations. \\

\section {Introduction.}
In the theory of differential equations for describing the behavior of the phase curves in a neighborhood of a singular point we use Lyapunov theory [1]. In the second method of Lyapunov  Lyapunov functions play a significant role. In [2], it was proposed an approach, based on the consideration of behavior of phase curves not along all the hypersurfaces of Lyapunov functions, but only for some of their sequences. In this paper we continue research in this direction.

\section { Embedded hypersurfaces }
\begin{defn}
Let ${\textbf{H}}^{n-1}$ be a coherent compact hypersurface (smooth compact $(n-1)$-dimensional submanifold ), which lies in a neighborhood $ G \subset \textbf R^n $, of the origin . We say that hypersurface $\textbf H ^{n-1}$ \textbf {restricts the origin}, if ${\textbf 0} \notin {\textbf H}^{n-1} $ and  ${\textbf H}^{n-1}$ is the bounded set of a compact set $ \textbf K $, containing the origin.
\end{defn}

\begin{defn}
Let $\textbf{H}_i^{n-1} \subset {\textbf R^n}$  be hypersurfaces, restricting the origin $ (i = 1, 2, ...) $. We say that $ \textbf H_i^{n-1} $ is a {\bfseries convergent sequence of hypersurfaces} if $\textbf H_i^{n-1}$ does not intersect and there are compact sets $ \textbf K_i $ from the outside $\textbf{H}_i^{n-1}$ for which the conditions: $ \bigcap_i \textbf{K}_i = \textbf{0} \in R ^ n$ and $ \textbf{K}_i \ supset \textbf{K}_j $ for $ i <j $.
\end{defn}

\begin{lemma}
Let  $ G $ be a neighborhood of the origin define in $ G $ a positive (negative)
differentiable function  $ z = F (\textbf{x}) $ such that $ F (\textbf{0}) = 0 $.
Then there is a convergent sequence of hypersurfaces $ \textbf{H} _i ^{n-1} $ in $ G $.
The numbers $ a_i = F (\textbf{H} _i ^{n-1}) $ are regular values of $ z = F (\textbf{x}) $. The sequence $ (a_i) $ converges  to $ 0 $ when $ i $ tends to $ \infty $.
\end{lemma}

{\it Proof.} Suppose that the function $ z = F (\textbf{x}) $ is  positivly defined in $ G $ . By theorem
of Sardis [3],we can choose a sequence $ (a_i) $ different regular values of 
$ Z = F (\textbf{x})s$, which correspond to zero monotonnically.
Consider the ball $ D_{\varepsilon} $ of radius $ \varepsilon $
centered at the origin. Let fix index $ i_o $ and construct the set
$ F ^{-1} [0, a_{i_o}] $, which may consist of many connected components. Select the connected component $ P_{i_o} $ of the origin. The set $ P_{i_o} $ is a smooth manifold with boundary, and in general,it may not be compact. We will show that there are functions $ a_{i_k}$, for which the corresponding connected component $ P_{i_k} $, of the origin lies inside the ball $ D_{\varepsilon} $.
Suppose by conradictions all values $ a_i $ corresponding connected components $ P_{i} $, which includes the origin, do not lie inside the ball $ D_{\varepsilon} $ and then crossing its border - the sphere $ S ^{ n-1} _{\varepsilon} $. Note that  for the construction for each index $ i $ the set $ Q_{i} = P_{i} \cap D_{\varepsilon} $ is not empty.
Consider contraction of functions $ z = F (\textbf{x}) $ on the sphere $ S ^{n-1} _{\varepsilon} $ and denote the function obtained by $ \widehat{F} $.
Obviously, $\widehat{F}$ is a differential function.
Since the sphere $S^{n-1} _{\varepsilon} $ is a compact smooth submanifold, then the value of the function
$ \widehat{F} $ is the segment $ [d_o, d_1] $,that $ 0 \notin [d_o, d_1] $.
Let $ a_{j_o} $ - be a value of $ z = F (\textbf{x}) $
strictly less than $ [d_0] $. Clearly, connected component  $ P_{j_o} $ lies inside the sphere $ S ^{n-1} _{\varepsilon} $, and it dose not cross it, because on it function $ Z = F (x_1, x_2, ..., x_n) $ takes values in the interval
 $ [0, a_{j_o}] $ and the origin belongs to the $ P_{j_o} $.
Thus edge manifold $ P_{j_o} $, which denotes by $ \textbf{H}_{i_0}^{n-1} $ is hypersurface that limits the origin. Clearly, for all other values of $ (a_j) $ functions $ z = F (\textbf{x}) $, which are less than $ a_{j_o} $ relevant connected components  $ P_{j} $ will be satisfy the condition $ P_{ j_o} \supset P_{j_o + 1} \supset P_{j_o + 2} \supset ... $.
By constructing limits manifolds $ P_{j} $ is a convergent sequence of hypersurfaces.
 For negatively we can argue likewise. $ \square $
\begin{defn}
Assume in the neighborhood $ G $, of the origin $ G $ is given a differentiable
 function $ z = F (\textbf{x}) $ such that
$ F (\textbf{0}) = 0 $.
Say that the function $ z = F (\textbf{x}) $ satisfies condition $ \textbf{L} $, if there is a sequence of regular values $ (a_i) $, which coincides with $ 0 $, set $ F ^{- 1} (a_i) $ has a connected component of $ \textbf{H}_{i} ^{n-1} $, which is a smooth hypersurface, that limits the origin and diameter of hypersurfaces $\textbf{H}_{i} ^{n-1} $ tends to $ 0 $ when $ i $ tends to infinity.
 \end{defn}
Since Lemma 2.1 implies that positive (negative) defined differentiable function that is specified in the neighborhood of origin, satisfies the condition $ \textbf{L}$.
Suppose that in  $ G $ is given a differentiable
function $ z = F(\textbf{x}) $, which satisfies the condition $ \textbf{L} $. Select a certain convergent sequence of hypersurfaces
$\textbf{H}_{i}^{n-1} $. At each point $\textbf{x} \in \textbf{H}_{i_0}^{n-1} $ define a nonzero gradient vector $ \overrightarrow{grad} F (\textbf{x}) $. We assign $ + $ for hypersurface if $ \overrightarrow{grad} F (\textbf{x}) $ is directed inside $\textbf{H}_{i_0} ^{n-1} $ and assign $ - $ otherwise.  The sign of hypersurface $\textbf{H}_{i_0} ^{n-1} $ will denote by $ \varepsilon ((\textbf{H}) _{i_0} ^{n-1}) $.
\begin{lemma}
Assume in the neighborhood of the origin  $ G $ is given a
differentiable function  $ z = F (\textbf{x}) $ such that $ F (\textbf{0}) = 0 $. Suppose that the origin is an isolated connected component of level surface $ F ^{-1} (0) $. Then the function $ z = F (\textbf{x}) $ satisfies condition $ \textbf{L}$.
\end{lemma}
{\it Proof.} Consider the ball $ D ^ n $ centered at the origin such that $ D ^ n \bigcap F ^{-1} (0) = \textbf{0}$.
Due to the isolation of the origin, connected components of the level surface $ F ^{-1} (0) $, so there exists always a ball.
Select the arbitrary point $\textbf{x}$ inside of $ D ^ n $ and let $ a = F (\textbf{x}) $. If $ a> 0 $, then for all points of the ball $ D ^ n $, except the origin, the function $ z = F (\textbf{x}) $ is positive. By contradiction , suppose there exists in  $D^n$ a point $\textbf {y}$ such that  $b = F(\textbf {y})$ and $ b < 0 $. Consider a ball in $ D ^ n $ and continuous path $ \gamma(t)$, which connects the points $ \textbf{x}$ and $ \textbf{y}$ and that does not go through the origin contraction the functions $ z = F (\textbf{x}) $ on $ \gamma(t)$ is a continuous function $ g (\gamma(t)) $, which takes on the ends of $ \gamma(t) $ opposite value signs.
By the theorem about the intermediate value of continious function, $ g (\gamma(t)) $ should take value $ 0 $, but it is impossible for the construction. The obtained contradiction proves the fact that in the ball $ D ^ n $ function $ z = F (\textbf{x}) $ is positive definite and therefore satisfies the condition $ \textbf{L} $.
If $ a <0 $ the similar arguments. $ \square $ \\

\section { Discrete analogue of a theorem by Lyapunov stability.}

Set on the hypersurfaces $\textbf{H}_i^{n-1} $  unit normal vector field $ \vec{N}(\textbf{x})$, which is directed to the inside of $ K_i $.

\begin{theorem}
Assume in the neighborhood of the origin $ G $ is given an autonomous system of ordinary differential equations
$ D (\textbf{x}) / dt = \textbf{f}(\textbf{x})$ $ (2.1) $, in which the origin is an isolated  equilibrium.
Suppose that in $ G $, there is a convergent sequence of hypersurfaces
$ \textbf{H}_i^{n-1} $.
If all points of $\textbf{x} \in \textbf{H}_i ^{n-1} $ a function $ S (\textbf{x}) = <\vec{N} (\textbf{x}) , \vec{f} (\textbf{x})> $ $(*)$ be positive, then the equilibrium of $ (2.1) $ is
 stable by Lyapunov. (We checked via $ \vec{f} (\textbf{x}) = (f_1 (\textbf{x}), ..., f_n (\textbf{x})).)$
\end{theorem}
{\it Proof.} It is obvious that assumptions of the theorem guarantee that phase curve $\textbf{X}(t) $ of $ (2.1) $, which begins at the point $\textbf{x}$
of hypersurface $\textbf{H}_{i_o} ^{n-1} $ is increasing at time $ t $ lies inside of the manifold $ K_{i_o} $, which serves as the border of $\textbf{H}_{i_o}^{n-1}$.
Otherwise, there exits this phase of the curve on manifold $ K_{i_o} $ requires that it crosses $\textbf{H}_{i_o} ^{n-1} $ along vector, which aims to look towards $\textbf{H}_{i_o}^{n-1} $, or is tangent to $\textbf{H}_{i_o}^{n-1}$.
But condition $(*)$ allows such behavior of the phase curve. It means stability and equilibrium
by Lyapunov. $ \square $

\begin{cor}
suppose in the $ G $ is given the system $ (2.1) $.
If the neighborhood of $ V \subset G $ is a differentiable function of $ z = F (\textbf{x}) $, which satisfies
the condition $ (\textbf{L}) $ and such that its derivative sign $ dF ((\textbf{x}) (t)) / dt $
along an arbitrary phase trajectory $\textbf{X}(t) $ of the system $(2.1)$ at points of hypersurfaces
$\textbf{H}_i^{n-1} $ coincides with the sign $\varepsilon(\textbf{H}_i^{n-1})$, then equilibrium
of system $ (2.1) $ is stable by Lyapunov.
\end{cor}

\section {Conditions of asymptotic stability by \\Lyapunov.}

\begin{defn}
Assume  in the  neighborhood $ G $, of the origin is fixed a convergent sequence of hypersurfaces
$\textbf{H}_{i}^{n-1} $. We say that hypersurfaces in sequence $\textbf{H}_{i} ^{n-1} $ {\bfseries are different}, if 
for every natural number $ n \in \textbf{N}$ there are numbers $ k, l \in \textbf{N}$ such that $ k> l> n $ hypersurfaces
 $\textbf{H}_{k} ^{n-1} $ and $\textbf{H}_{l} ^{n-1} $ are not homeomorphic.
\end{defn}

\begin{theorem}
Let in the  $ G $ is given the system $ (2.1) $ and there is a convergent sequence of
hypersurfaces  $\textbf{H}_i ^{n-1} $, which is different.
If all values $\textbf{x} \in \textbf{H}_i ^{n-1} $, of the function
 $ S (\textbf{x}) = <\vec{N} (\textbf{x}), \vec{f} (\textbf{x})> $ are positive, then the equilibrium of $ (01.02) $
is stable but not asymptotically stable by Lyapunov.
\end{theorem}

{\it Proof.}
 By theorem 3.1 the origin is stable according to Lyapunov. Suppose that the origin 
is asymptotically stable by Lyapunov.Then each phase curve  of $ (2.1) $, which begins 
at a fixed hypersurface $\textbf{H}_{i_0}^{n-1} $ with an increase in the parameter $ t $ crossing 
each hypersurface  $\textbf{H}_{j} ^{n-1} $ $ (j> i_0) $ only at one point. That phase 
curve of $ (2.1) $ defines a homeomorphism between hypersurface $\textbf{H}_{i} ^{n-1} $.
But this is impossible, because hypersurface in the sequence $\textbf{H}_{i} ^{n-1} $ are different 
The obtained contradiction proves the theorem.
$\square$

\begin{rem}
Using Homologous not simply connected sphere $ \Sigma ^{n-1} $ $ (n> 3) $ one can construct convergent sequence of hypersurfaces $\textbf{H}_i ^{n-1} $ in $ G $, which is different. Among all $\textbf{H}_i ^{n-1} $ will be hypersurfaces that homeomorphic both the standard $ S ^{n-1} $ and homologous spheres $ \Sigma^{n-1} $. For such hypersurfaces $\textbf{H}_i ^{n-1} $ it is easy to set in $ G $ the system $(2.1) $, which satisfies the conditions of theorem 4.1.
\end{rem}

\begin{theorem}
Let in $ G $ there is the system $ (2.1) $ and there is a convergent sequence of
hypersurfaces $\textbf{H}_i ^{n-1}$. If all points  $\textbf{x} \in \textbf{H}_i ^{n-1} $  of the function
$ S(\textbf{x}) = <\vec{N} (\textbf{x}), \vec{f} (\textbf{x})> $ is positive and the origin is the only invariant set, then the  equilibrium of $ (2.1) $ is asymptotically stable by Lyapunov.
\end{theorem}

{\it Proof.} According to theorem 3.1 the origin is stable by Lyapunov. For each phase trajectory $ \gamma(t) $,
passing in the  neighborhood of the origin of $ \omega (\gamma(t)) $-limit sets are invariant sets.
So $ \gamma(t) $ approches to the origin, when $ t \rightarrow \infty $. The obtained contradiction proves the theorem. $ \square $

1. ~ {\it ~  N.Rush, P. Abets, M. Lalua.} Direct Lyapunov method in stability theory. Mir Pub,1980. p 300.

2. ~ {\it ~Yu. Sharko} Discrete conditions for stability by Lyapunov Proceedings of Institute of Mathematics, National Academy of Sciences of Ukraine. 2005,Vol 2,№ 3. p. 279-288.

3. ~ {\it ~ M. Hirsch } Differential topology Mir Pub,1979.  p.279.
\\\\
Іnstitute of mathematics of Ukrainan academy of sciences
\\
\end{document}